\documentclass{amsart}
\usepackage{amsmath, amssymb,color}
\newtheorem{Theorem}{Theorem}[section]

\newtheorem{Proposition}[Theorem]{Proposition}
\theoremstyle{definition}
\newtheorem{Definition}[Theorem]{Definition}
\theoremstyle{remark}
\newtheorem{rem}[Theorem]{Remark}
\numberwithin{equation}{section}

\newcommand{\R}{\mathbb R}
\newcommand{\N}{\mathbb N}

\newcommand{\C}{\mathcal{C}}

\newcommand{\D}{\mathcal{D}}
\newcommand{\F}{\mathcal{F}}
\newcommand{\p}{\mathcal{P}}
\newcommand{\A}{\mathcal{A}}
\newcommand{\id}{Id_{]0;1[}}
\begin{document}

\title{The group of diffeomorphisms of a non compact manifold is not regular.}%
\author{Jean-Pierre Magnot}

\address{LAREMA - UMR CNRS 6093 \\ Universit\'e d'Angers \\ 2 Boulevard Lavoisier -
	49045 Angers cedex 01 \\ and
	Lyc\'ee Jeanne d'Arc \\
         Avenue de grande bretagne \\
         F-63000 Clermont-Ferrand}
\email{jean-pierr.magnot@ac-clermont.fr}

\begin{abstract}
We show that a group of diffeomorphisms $\D$ on the open unit interval $I,$ equipped with the topology of uniform convergence on any compact set of the derivatives at any order, is non regular: the exponential map is not defined for some path of the Lie algebra. this result extends to the group of diffeomorphisms of finite dimensional, non compact manifold $M.$   
\end{abstract}

\maketitle



\section*{Introduction}
In the theory of compact Lie groups, the exponential map defines a chart around the neutral element. This map plays a central role in many aspects of the theory. In infinite dimensional Lie groups, the existence of the exponential map is not straightforward, mostly because of the lack of compact neighbourhood of the neutral element. {
	For this reason} we often work, in infinite dimensional Lie groups, only with \textbf{regular} Lie groups, which are Lie groups that have an exponential map, which is a map which integrates any smooth path $v$ on the Lie algebra
to a smooth path $g$ on the Lie group via the equation on  logarithmic derivatives $$dg.g^{-1} = v.$$ In easy examples, such as Banach Lie groups,
the exponential map also defines a chart around the neutral element, but in some other examples such as groups of diffeomorphisms on a compact manifold, it is not the case. This last pathology generates technical difficulties 
{
\cite{Om1973,Om1981}.}
 
{
	The notion of regular Lie group was first described by Omori, motivated by the discovery of non-enlargeable Lie algebras of vector fields \cite{Om1981} after a serie of works on the ILH structures of the group of diffeomorphisms of a compact manifold( see e.g. \cite{Om1973})} The terminology of ``regular Lie group'' was introduced with a second class of examples: groups of Fourier-integral operators {
	in a serie of papers . For an organized exposition, see \cite{Om}.}
Since these founding examples, despite of many efforts, {
	there is still no known example of non-regular Fr\'echet Lie group}. As quoted in \cite{KM}, there exists many candidates, but it is quite uneasy to show that a differential equation on an infinite dimensional Lie group has no solution. 

{
	In addition, we have to precise that the universal setting for infinite dimensional geometry seems not to exist. Historically, various authors developped more and more general settings for infiniet dimensional ``manifolds'': Hilbert, Banach, then Fr\'echet and locally convex manifolds (equipped with atlas), and then raised the necessity to deal with ``manifolds'' without charts. This phenomenon was even so embarrassing that the precise setting for differential calculus on infinite dimensional setting has been skipped by several authors for applications, see e.g. \cite{Khe}. Several definitions and settings have been given by various authors, and the choice that we make to use diffeological spaces as a maximal category, and Fr\"olicher spaces as an intermediate category for differential geometry \cite{Ma2006-3,Wa}, is still quite controversial but this setting becomes developped enough to raise applications. Moreover, historically, diffeological spaces have been developped by Souriau in the 80's with the motivation to deal with the precise objects of interest here, that is groups of diffeomorphisms on non compact, locally compact, smooth manifolds without boundary.} 
 
In this short communication, we describe an example of non regular group of diffeomorphisms 
{
	$\D$} on the unit interval, which is {
	in the setting of}  Fr\"olicher Lie groups. The difference is that, in the Fr\"olicher setting, existence of charts is not assumed. For basics on this setting, due to Fr\"olicher and Kriegl, see \cite{FK,KM} and e.g. \cite{Ma2006-3,Ma2013,Wa} for a short exposition on Fr\"olicher spaces and Fr\"olicher Lie groups. 

This example, which can appear as a toy example, {
	acts as a preliminary result in order to show that the (full) group of diffeomorphisms $Diff(M)$ of a non-compact, finite dimensional manifold $M$, equipped with the topology of uniform convergence of any partial derivative on any compact subset of $M$, called the smooth compact-open topology (or weak topology in \cite{Hir}), is non regular. This topology appears as quite natural, and quite minimal compared to the (much stronger) topologies described in \cite[section 43.1]{KM} (similar to the classical $C^\infty$-Whitney topology) and in the more diversified \cite{KMR}, which furnish examples of regular Lie groups (with atlas). In our example, no chart can be actually successfully constructed by the lack of adequate implicit functions theorem in this topology, and also because this example seems not locally homeomorphic to any open subset of a function space for the considered (smooth compact-open) topology. Of course, these considerations are still open questions, since from another viewpoint no-one can prove actually neither that this is impossible, nor that it is possible, to get these local homeomorphism properties for $\D.$

Our method of proof} is inspired by:

- first the example of non integrable Lie subalgebra of $Vect(S^1)$ due to Omori {
	\cite{Om1981},} where non integrability is due to the existence of translations,

- secondly the non existence of horizontal lifts for connections on fiber bundles when the typical fiber is not compact, see \cite{KMS}. Moreover, if the typical fiber is compact, connections on the considered fiber bundle are in one-to-one correspondence with connections on a principal bundle with structure group a (regular) group of diffeomorphisms, see e.g. \cite{KM}.

We need to precise that the example that we develop seems already known in the mathematical literature,{
	but nowhere stated clearly to our knowledge so} that we feel the need of a rigorous description of the announced phenomenon: a constant path on the Lie algebra is not integrable into a path of the group.
This is done in three steps:
\begin{itemize}
\item first summarize the settings that have been developed to enable a rigorous differential geometry on groups of diffeomorphisms on non compact manifolds (section 1), namely diffeological paces and Fr\"olicher spaces, applied to our example
of group of diffeomorphisms on the open unit interval, 
\item secondly show that the constant map on the unit interval is in {
	$T_{\id}\D$, and that this element understood as a constant path cannot integrate in $\D$ by an argument of ``translation semi-group'' discovered by \cite{Om1981} for a different goal.}
\item finally {
	embedding $\D$ into $Diff(M),$ for a non compact manifold $M,$ in a way such that the translation semi-group required on $\D$ cannot be described as a semi-group of transformations on $M,$ we get a direct proof, new to our knowledge, of the following result :}   
\end{itemize}

\begin{Theorem}
	Let $M$ be a locally compact, non compact manifold. Then the group $Diff(M),$ equipped with its {
		functional} diffeology, is a diffeological Lie group which is non regular. 
\end{Theorem}

\section{Preliminaries}
\subsection{Souriau's diffeological spaces and Fr\"olicher spaces}
\label{1.1}

\begin{Definition} \cite{Sou}, see e.g. \cite{Igdiff}. Let $X$ be a set.
	
	\noindent $\bullet$ A \textbf{p-parametrization} of dimension $p$ (or $p$-plot)
	on $X$ is a map from an open subset $O$ of $\R^{p}$ to $X$.
	
	\noindent $\bullet$ A \textbf{diffeology} on $X$ is a set $\p$
	of parametrizations on $X$, called plots of the diffeology, such that, for all $p\in\N$,
	
	- any constant map $\R^{p}\rightarrow X$ is in $\p$;
	
	- Let $I$ be an arbitrary set of indexes; let $\{f_{i}:O_{i}\rightarrow X\}_{i\in I}$
	be a family of compatible maps that extend to a map $f:\bigcup_{i\in I}O_{i}\rightarrow X$.
	If $\{f_{i}:O_{i}\rightarrow X\}_{i\in I}\subset\p$, then $f\in\p$.
	
	- Let $f\in\p$, defined on $O\subset\R^{p}$. Let $q\in\N$,
	$O'$ an open subset of $\R^{q}$ and $g$ a smooth map (in the usual
	sense) from $O'$ to $O$. Then, $f\circ g\in\p$.
	
	\vskip 6pt $\bullet$ If $\p$ is a diffeology on $X$, then $(X,\p)$ is
	called a \textbf{diffeological space}.
	
	\noindent Let $(X,\p)$ and $(X',\p')$ be two diffeological spaces;
	a map $f:X\rightarrow X'$ is \textbf{differentiable} (=smooth) if
	and only if $f\circ\p\subset\p'$. \end{Definition}

\begin{rem}
	Any diffeological space $(X,\p)$ can be
	endowed with the  weakest topology such that all the maps that
	belong to $\p$ are continuous. {
		This topology is called D-topology, see \cite{CSW}.}
\end{rem}

We now introduce Fr\"olicher spaces, see \cite{FK}, using the
terminology defined in \cite{KM}.

\begin{Definition} $\bullet$ A \textbf{Fr\"olicher} space is a triple
	$(X,\F,\C)$ such that
	
	- $\C$ is a set of paths $\R\rightarrow X$,
	
	- $\F$ is the set of functions from $X$ to $\R$, such that a function
	$f:X\rightarrow\R$ is in $\F$ if and only if for any
	$c\in\C$, $f\circ c\in C^{\infty}(\R,\R)$;
	
	- A path $c:\R\rightarrow X$ is in $\C$ (i.e. is a \textbf{contour})
	if and only if for any $f\in\F$, $f\circ c\in C^{\infty}(\R,\R)$.
	
	\vskip 5pt $\bullet$ Let $(X,\F,\C)$ and $(X',\F',\C ')$ be two
	Fr\"olicher spaces; a map $f:X\rightarrow X'$ is \textbf{differentiable}
	(=smooth) if and only if $\F'\circ f\circ\C\subset C^{\infty}(\R,\R)$.
\end{Definition}

Any family of maps $\F_{g}$ from $X$ to $\R$ generates a Fr\"olicher
structure $(X,\F,\C)$ by setting, after \cite{KM}:

- $\C=\{c:\R\rightarrow X\hbox{ such that }\F_{g}\circ c\subset C^{\infty}(\R,\R)\}$

- $\F=\{f:X\rightarrow\R\hbox{ such that }f\circ\C\subset C^{\infty}(\R,\R)\}.$

In this case we call $\F_g$ a \textbf{generating set of functions}
for the Fr\"olicher structure $(X,\F,\C)$. One easily see that
$\F_{g}\subset\F$. This notion will be useful in the sequel to
describe in a simple way a Fr\"olicher structure, see for instance
Proposition \ref{froproj} below. A Fr\"olicher space $(X,\F,\C)$
carries a natural topology, which is the pull-back topology of
$\R$ via $\F$. We note that in the case of a finite dimensional
differentiable manifold $X$ we can take $\F$ the set of all smooth
maps from $X$ to $\R$, and $\C$ the set of all smooth paths from
$\R$ to $X.$ In this case the underlying topology of the
Fr\"olicher structure is the same as the manifold topology
\cite{KM}. In the infinite dimensional case, there is to our
knowledge no complete study of the relation between the
Fr\"olicher topology and the manifold topology; our intuition is
that these two topologies can differ.

We also remark that if $(X,\F, \C)$ is a Fr\"olicher space, we can
define a natural diffeology on $X$ by using the following family
of maps $f$ defined on open domains $D(f)$ of Euclidean spaces
(see \cite{Ma2006-3}):
$$
\p_\infty(\F)=
\coprod_{p\in\N}\{\, f: D(f) \rightarrow X; \, \F \circ f \in C^\infty(D(f),\R) \quad \hbox{(in
	the usual sense)}\}.$$

If $X$ is a differentiable manifold, this diffeology has been
called the { \bf n\'ebuleuse diffeology} by J.-M. Souriau, see
\cite{Sou},{
	or nebulae diffeology in \cite{Igdiff}.} We can easily show the following:

\begin{Proposition} \label{fd} \cite{Ma2006-3}
	Let$(X,\F,\C)$
	and $(X',\F',\C')$ be two Fr\"olicher spaces. A map $f:X\rightarrow X'$
	is smooth in the sense of Fr\"olicher if and only if it is smooth for
	the underlying diffeologies $\p_\infty(\F)$ and $\p_\infty(\F').$
\end{Proposition}

Thus, we can also state:
\vskip 12pt

\begin{tabular}{ccccc}
	smooth manifold  & $\Rightarrow$  & Fr\"olicher space  & $\Rightarrow$  & Diffeological space
\end{tabular}

\vskip 12pt A deeper analysis of these implications has been given
in \cite{Wa}. The next remark is inspired on this work and on
\cite{Ma2006-3}; it is based on \cite[p.26, Boman's theorem]{KM}.
\begin{rem}
	We notice that the set of contours $\C$ of the Fr\"olicher space
	$(X,\F,\C)$ \textbf{does not} give us a diffeology, because a diffelogy
	needs to be stable under restriction of domains. {
		In the case of paths in
	$\C$ the domain is always $\R$ where as the domain of 1-plots can (and has to) be any interval of $\R.$} However, $\C$ defines a ``minimal diffeology''
	$\p_1(\F)$ whose plots are smooth parameterizations which are locally of the
	type $c \circ g,$ where $g \in \p_\infty(\R)$  and $c \in \C.$ Within this setting,
	we can replace $\p_\infty$ by $\p_1$ in Proposition \ref{fd}.
\end{rem}

We also remark that given an algebraic structure, we can define a
corresponding compatible diffeological structure. For example 
{
following \cite[p.66-68]{Igdiff}, a
$\R-$vector space equipped with a diffeology is called a
diffeological vector space if addition and scalar multiplication
are smooth.} An
analogous definition holds for Fr\"olicher vector spaces. 

\begin{rem} \label{comp}
{
	Fr\"olicher and Gateaux smoothness are the same notion
	if we restrict to a Fr\'echet context}.
	Indeed, for a smooth map $f : (F, \p_1(F)) \rightarrow \R$ defined
	on a Fr\'echet space with its 1-dimensional diffeology, we have
	that $\forall (x,h) \in F^2,$ the map $t \mapsto f(x + th)$ is
	smooth as a classical map in $\C^\infty(\R,\R).$ And hence, it is
	Gateaux smooth. The converse is obvious.
\end{rem}









\vskip 12pt
{
Diffeologies on cartesian products, projective limits, quotients, subsets as well as pull-back  and push-forward diffeologies are described in \cite[Chapter 1]{Igdiff}. The reader can refer also to \cite{Sou} or \cite{Ma2013,Ma2015} for faster exposition.}
\subsection{{
		Functional} diffeology}
Let $(X,\p)$ and $(X',\p')$ be two diffeological spaces. Let $M \subset C^\infty(X,X')$ be a set of smooth maps. The \textbf{{
		functional} diffeology} on $S$ is the diffeology $\p_S$
made of plots
$$ \rho : D(\rho) \subset \R^k \rightarrow S$$
such that, 
for each $p \in \p, $
the maps $\Phi_{\rho, p}: (x,y) \in D(p)\times D(\rho) \mapsto \rho(y)(x) \in X'$ are plots of $\p'.$ With this definition, we have the classical fundamental properties:

\begin{Proposition} \cite{Igdiff}
	Let $X,Y,Z$ be diffeological spaces, 
	$$C^\infty(X\times Y,Z) = C^\infty(X,C^\infty(Y,Z)) = C^\infty(Y,C^\infty(X,Z))$$
	as diffeological spaces equipped with {
		functional} diffeologies.
\end{Proposition}

\subsection{Tangent space}

There are actually two main definitions, (2) and (3) below, of the tangent space of a diffeological space:

\begin{enumerate}
	\item the \textbf{internal tangent cone} \cite{DN2007-1}. For each $x\in X,$ we consider $$C_{x}=\{c \in C^\infty(\R,X)| c(0) = x\}$$ and take the equivalence relation $\mathcal{R}$ given by $$c\mathcal{R}c' \Leftrightarrow \forall f \in C^\infty(X,\R), \partial_t(f \circ c)|_{t = 0} = \partial_t(f \circ c')|_{t = 0}.$$
	The internal tangent cone at $x$ is the quotient $$^iT_xX = C_x / \mathcal{R}.$$ If $X = \partial_tc(t)|_{t=0} \in {}^iT_X, $ we define the simplified notation  $$Df(X) = \partial_t(f \circ c)|_{t = 0}.$$
	\item The \textbf{internal tangent space} at $x \in X$ {
		described in}  \cite{CW}
	\item the \textbf{external} tangent space $^eTX,$ defined as the set of derivations on $C^\infty(X,\R).$ \cite{KM,Igdiff}.
\end{enumerate}

It is shown in  \cite{DN2007-1} that the internal tangent cone at a point $x$ is not a vector space in many examples. This motivates  \cite{CW}. For finite dimensional manifold, definitions (1), (2) and (3) coincide. {
	For more comparisons, see \cite[section 28]{KM} for a comparison for infinite dimensional manifolds and also \cite{CW}.}
\subsection{Regular Lie groups}

\begin{Definition}
	Let $G$ be a group, equiped with a diffeology $\p.$ We call it \textbf{diffeological group} if both multiplication and inversion are smooth.
\end{Definition}

The same definitions hold for Fr\"olicher groups.
{
	Let us now recall \cite[Proposition
	1.6.]{Les}, which shows that the distinction between internal tangent cone and internal tangent space is not necessary for diffeological groups.
	\begin{Proposition} \label{leslie}
		Let $G$ be a diffeological group. Then the tangent cone at the neutral element $T_eG$ is a diffeological vector space.
\end{Proposition}}
Following Iglesias-Zemmour, \cite{Igdiff}, who does not assert that arbitrary diffeological groups have a Lie algebra,
{
we restrict ourselves to a smaller class of diffeological groups } which have such a tangent space at the neutral element.
{
	Intuitively speaking, The diffeological group $G$ is a \textbf{diffeological Lie group} if and only if
	the derivative of the Adjoint action of $G$ on $^iT_eG$ defines a Lie bracket. In this case, we call $^iT_eG$ the Lie algebra of $G$, that we note generically $\mathfrak{g}.$ One crucial question consists in giving a technical condition which ensures the classical properties of Adjoint and adjoint actions, e.g.:
\begin{itemize}
	\item Let $(X,Y) \in \mathfrak{g}^2,$ $X+Y = \partial_t(c.d)(0)$  where $c,d \in \C ^2,$ $c(0) = d(0) =e_G ,$
	$X = \partial_t c(0)$ and $Y = \partial_t d(0).$
	\item Let $(X,g) \in \mathfrak{g}\times G,$ $Ad_g(X) = \partial_t(g c g^{-1})(0)$  where $c \in \C ,$ $c(0) =e_G ,$
	and $X = \partial_t c(0).$
	\item Let $(X,Y) \in \mathfrak{g}^2,$ $[X,Y] = \partial_t( Ad_{c(t)}Y)$   where $c \in \C ,$ $c(0) =e_G ,$
	$X = \partial_t c(0).$
\end{itemize}

According to \cite{Ma2013} (which deals with Fr\"olicher Lie groups), one can assume only that the desired properties are fulfilled, leaving technicities for specific examples. One criteria has been given in \cite[definition 1.13 and Theorem 1.14]{Les} but this is not necessary here since in the framework that we consider, the properties of the Lie bracket will rise naturally and directly.
For these reasons, we
 }
give the following definition:

\begin{Definition}
	The diffeological group $G$ is a \textbf{diffeological Lie group} if and only if
	the derivative of the Adjoint action of $G$ on $^iT_eG$ defines a smooth Lie bracket. In this case, we call $^iT_eG$ the Lie algebra of $G$, that we note generically $\mathfrak{g}.$
\end{Definition}

Let us now concentrate on diffeological Lie groups, and in this case we note $\mathfrak{g}={}^iT_eG.$
The basic properties of adjoint, coadjoint actions, and of Lie brackets, remain globally the same
as in the case of finite-dimensional Lie groups, and the proofs are similar: we
only need to replace charts by plots of the underlying diffeologies (see e.g. \cite{Les} for further details, and \cite{BN2005} for the case of Fr\"olicher Lie groups), as soon as one has checked that the Lie algebra $\mathfrak{g}$ is a diffeological Lie algebra, i.e. a diffeological vector space with smooth Lie bracket.

\begin{Definition} \label{reg1} \cite{Les} A diffeological Lie group $G$ with Lie algebra $\mathfrak{g}$
	is called \textbf{regular} if and only if there is a smooth map \[
	Exp:C^{\infty}([0;1],\mathfrak{g})\rightarrow C^{\infty}([0,1],G)\]
	such that $g(t)=Exp(v(t))$ is the unique solution
	of the differential equation \begin{equation}
	\label{loga}
	\left\{ \begin{array}{l}
	g(0)=e\\
	\frac{dg(t)}{dt}g(t)^{-1}=v(t)\end{array}\right.\end{equation}
	We define the exponential function as follows:
	\begin{eqnarray*}
		exp:\mathfrak{g} & \rightarrow & G\\
		v & \mapsto & exp(v)=g(1) \; ,
	\end{eqnarray*}
	where $g$ is the image by $Exp$ of the constant path $v.$ \end{Definition}

	



\subsection{Groups of diffeomorphisms}
Let $M$ be a localy compact, non compact manifold,
{
which is assumed Riemannian without restriction, 
} equipped with its n\'ebuleuse diffeology. We equip the groups of diffeomorphisms $Diff(M)$ with the topology of convergence of the derivatives an any order, uniformly on each compact subset of $M,$ usually called $C^\infty-$compact-open topology or weak topology in \cite{Hir}. Traditionnally, $Vect(M)$ is given as the Lie algebra of $Diff(M),$ but \cite[section 43.1]{KM} shows that this strongly depends on the topology of $Diff(M).$ 
{
Indeed, the Lie algebra of vector fields described in \cite[section 43.1]{KM} is the Lie algebra of compactly supported vector fields, which is not the (full) Lie algebra $Vect(M).$
 In another context, when $M$ is compact, $Vect(M)$ is the Lie algebra of $Diff(M),$ which can be obtained by Omori's regularity theorems \cite{Om1973,Om} and recovered in \cite{CW}. 
}
What is well known is that infinitesimal actions 
{
(i.e. elements of the internal tangent space at identity)
} of $Diff(M)$ on $C^\infty(M,\R)$ generates vector fields, viewed as order 1 differential operators.
{
The bracket on vector fields is given by $$(X,Y)\in Vect(M)	\mapsto [X,Y]= \nabla_XY - \nabla_YX$$
where $\nabla$ is the Levi-Civita connection on $TM.$ This is a Lie bracket, stable under the Adjoint action of $Diff(M).$
} Moreover, the compact-open topology on $Diff(M)$ generates a corresponding $C^\infty-$compact-open topology on $Vect(M).$ This topology is itself the $D-$topology for the the {
functional} diffeology on $Diff(M).$  Following \cite[Definition 1.13 and Theorem 1.14]{Les}, $Vect(M)$ equipped with the $C^\infty$ compact-open topology is a Fr\'échet vector space, {
and  the Lie bracket is smooth.
Moreover we feel the need to remark that}  the evaluation maps
$$ T^*M \times Vect(M) \rightarrow \R$$
separate $Vect(M).$ Thus $Diff(M)$ is a diffeological Lie group 
{
matching with the criteria of \cite[Definition 1.13 and Theorem 1.14]{Les}, and}
 for the 
{
functional 
}
diffeology, with Lie algebra
{
$\mathfrak{g} \subset Vect(M).$}
\section{A non regular group of diffeomorphims of the unit interval}
\subsection{Premilinaries}
Let $F$ be the vector space of smooth maps $f \in C^\infty(]0;1[;\R).$ We equip $F$ with the following semi-norms:

 For each $(n,k) \in \N^* \times \N $,  $$||f||_{n,k} = \sup_{\frac{1}{n+1}\leq x \leq \frac{n}{n+1}} |D^k_xf|.$$

This is a Fr\'echet space, and its topology is the smooth compact-open topology, which is the $D-$topology of the compact-open diffeology. 
Let 
$$\A = \{ f \in C^\infty(]0;1[;]0;1[)|\lim_{x \rightarrow 1}f(x)=1
 \wedge\lim_{x \rightarrow 0}f(x)=0 \}.$$
Finally, we set
 $$\D = \{ f \in \A | \inf_{x \in ]0;1[}f'(x) >0 \}.$$
 $\D$ is a {
 	contractible} set of diffeomorphisms of the open interval $]0;1[$ which is an (algebraic) group for composition of functions. Composition of maps, and inversion, is smooth for the {
 	functional} diffeology.
Unfortunately, $\D$ is not open in $\A.$
As a consequence, we are unable to prove that it is a Fr\'echet Lie group. However, considering the smooth diffeology induced on $\D$ by $\A,$ the inversion is smooth. As a consequence, $\D$ is (only) a diffeological Lie group.
  
\subsection{A non integrable path of the Lie algebra}
{
Let us consider the standard mollifier  $$\phi(u) = \frac{1}{K}e^{\frac{1}{u^2-1}}$$ defined for $ u \in ]-1;1[$, with $$K = \int_{-1}^{1}e^{\frac{1}{u^2-1}}du >0,4$$ and extended smoothly to $\R$ setting $\phi(u) = 0 $ whenever $u \notin ]-1;1[.$ We set $\phi_\alpha(u) = \frac{1}{\alpha}\phi\left(\frac{x}{\alpha}\right).$ Let us define, with $\mathbf{1}$ the standard characteristic function of a set,
\begin{eqnarray*}
	c_t(x) & = &x + t \left(\mathbf{1}_{[|t|; 1-|t|]}(x)*\phi_{|t|}(x)\right)\\ &= &\left\{\begin{array}{cl}\forall x<2|t|, & c_t = x + t\int_{|t|-x}^{|t|}\phi_{|t|}(u)du\\ \forall x \in [2|t|; 1-2|t|], & c_t(x) = x+t \\
		\forall x>1-2|t|, & c_t = x + t\int_{-|t|}^{1-|t|-x}\phi_{|t|}(u)du
	\end{array}\right.\end{eqnarray*}
where $*$ is the standard convolution in the $x-$variable, and $-1/4 < t< 1/4.$ For $t \neq 0,$ $(t,x) \mapsto c_t(x))$ is smooth.  Moreover, since $$\sup \phi < \frac{1}{0,4e} <1,$$ for fixed $t,$ one can easily check that $\partial_x c_t > 0,$ which shows that $\forall t \in ]-1/4;1/4[, $ $c_t \in \D.$

\begin{Theorem}\label{path}
	The path $ t \mapsto c_t $ is of class $C^\infty(]-1/4;1/4[,\D).$ Moreover, 
	$$ \partial_t c_t|_{t=0} = \mathbf{1}_{]0;1[}.$$
	is a constant map.
\end{Theorem}

\noindent
\textbf{Proof.} 
Let $x \in ]0;1[.$
Let us now check smoothness of $(t,x) \mapsto c_t(x)$ at $(0;x).$ Let $\alpha = \min\{x,1-x\}$ and let us restrict our study on the open subset $\left] -\frac{\alpha}{4};\frac{\alpha}{4}\right[\times\left]x - \frac{\alpha}{2} ; x+ \frac{\alpha}{2}\right[.$  For $t \in \left] -\frac{\alpha}{4};\frac{\alpha}{4}\right[,$ we have that $$c_t(u)= x+t$$ whenever $$u \in\left]x - \frac{\alpha}{2} ; x+ \frac{\alpha}{2}\right[\subset\left]\frac{\alpha}{2} ; 1- \frac{\alpha}{2}\right[\subset \left]2|t|;1-2|t|\right[.$$
By the way, $(u,t) \mapsto c_t(u)$ is smooth on $\left] -\frac{\alpha}{4};\frac{\alpha}{4}\right[\times\left]x - \frac{\alpha}{2} ; x+ \frac{\alpha}{2}\right[$
and hence it is smooth everywhere. By direct differentiation, 

$$\partial_tc_t (x)|_{t=0} = 1.$$\qed}

Now, we get the following:
\begin{Theorem} \label{nr}
There exists a smooth path $v $ on $T_{\id} \D$ such that no smooth path $g $ on $\D$ satisfies the equation $$\partial_t g \circ g^{-1} == v.$$
\end{Theorem}

\noindent
\textbf{Proof.}
Let $v$ be the constant path equal to $\mathbf{1}.$ Let $t \mapsto g_t$ be a solution of the last equation and let $x \in ]0;1[.$  
{
Then we have, $\forall   y \in ]0;1[,$ setting $x = g_t(y),$ $$\partial_t g_t (y)= \left(\partial_t g_t\right) \circ g_t^{-1}(x) = 1,$$ and by the way, $$\forall(t,y) \in \left] -1/4;1/4\right[ \times \left] 0;1 \right[, \partial_t c_t(y) = 1 \hbox{ and } c_o(y)=y.$$
} so that the only possible solution is  the translation $g_t(x) = x+t.$ We have $$\forall t > 0, g_t \notin \A$$ so that $$\forall t>0, g_t \notin \D.$$ 
\qed

As a consequence, we get the announced result:
\begin{Theorem}
	$\D$ is a non regular diffeological Lie group. 
	\end{Theorem}

\subsection{Final remark: $Diff(M)$ is a non regular diffeological Lie group}
{
}
 Let $M$ be a smooth manifold,  for which there is an embedding 
 $$ e : ]0;1[  \rightarrow E$$ for which $Im(e)$ is closed in $M.$
 We also assume that this embedding can be extended to a so-called ``thick path", which image will be a tubular neighbourhood of $Im(e)$. These are the main conditions that are needed, which are fulfilled when $M$ is $n-$dimensional. Under these conditions, following the "smooth tubular neighbourhood theorem" (see e.g. \cite{B}) we can assume that there is a parametrization of the closed tubular neighbourhood under consideration via an embedding
 $$ E:  ]0;1[\times B_{n-1} \rightarrow M$$
 where $B_{n-1}$ is the Euclidian $n-1$ dimensional unit ball and $E(x,0) = e(x).$ We set $phi: \R \rightarrow \R_+$ a smooth function with support in $[-1;1]$ and such that $\phi = 1$ on a neighbourhood of $0.$
 We also parametrize $B_{n-1}$ via spherical coordinates $(r, \theta_1,...\theta_{n-1})$
 Under these conditions, we consider the path
 $$ C_t : \R \rightarrow Diff(M)$$
 defined by 
 $$C_t(x)= \left\{\begin{array}{cl}  x & \hbox{if } x \notin Im(E) \\
 E(c_{\phi(r).t}(x'),r, \theta_1, ... \theta_{n-1}) & \hbox{if } x = E(x',r, \theta_1, ... \theta_{n-1})\end{array}\right. .$$ 
 We have that $ \partial_t C_t(x)|_{t=0}$ is a smooth vector field, which equals to $e_* (\mathbf{1}) $ when $r=0.$
 As a consequence,
 
 \begin{Theorem}
 	The vector field $ \partial_t C_t(x)|_{t=0}$ has no global flow on $M$ and hence $Diff(M)$ is a non regular diffeological Lie group with integral Lie algebra.
 \end{Theorem}
 
 \vskip 12pt
 \noindent
 \textbf{Proof.}
 We investigate the flow of $ \partial_t C_t(x)|_{t=0}$ on $Im(e)$
 and we have that $$e^*(\partial_t C_t(x)|_{t=0}) = \mathbf{1} \times \{0\} \in e^*TM = C^\infty(]0;1[; \R \times \R^{n-1}).$$ Thus the flow along $Im(e)$ must be $t \mapsto e(x+t)$. Since $Im(e)$ is closed in $M$; this flow does not extend to a flow in $M.$ \qed
 
 \vskip 12pt
 \noindent
 {
 	\textbf{Acknowledgements:} the author would like to thank 
 	Daniel Christensen, who helped to improve the last version of this paper with very interesting comments and questions.}

\end{document}